\pdfoutput=1 
\documentclass[11pt,reqno]{amsart}

\usepackage[utf8]{inputenc} 
\usepackage[noadjust]{cite}
\usepackage{amsmath,amssymb,amsthm}
\usepackage{upgreek}
\usepackage{graphicx}

\usepackage{booktabs}
\usepackage{siunitx} 
\usepackage{caption}
\usepackage{makecell} 
\usepackage{placeins} 

\usepackage{listings}
\usepackage{xcolor} 
\usepackage[T1]{fontenc}

\usepackage{lmodern}

\usepackage{microtype}

\usepackage{xurl}

\usepackage{hyperref}
\hypersetup{breaklinks=true}

\renewcommand{\tau}{\uptau}

\newcommand{\tautitle}{\texorpdfstring{\ensuremath{\boldsymbol{\tau}}}{tau}}

\begin{document}
\raggedbottom 

\title[Another Multivector Inverse Algorithm]{Another Multivector Inverse Algorithm}

\author[D.~A.~Jones]{David Arthur Jones}
\address{Independent researcher\\Springfield, VA, USA}
\email{dajonesy@gmail.com}
\urladdr[ORCID]{https://orcid.org/0009-0009-7693-3396}


\begin{abstract}
A numerical algorithm is described for computing the inverse of a non-singular multivector across dimensions.  The method wraps a reversion-like product around a Faddeev--LeVerrier--Souriau (FLS) polynomial core.  While FLS provides dimensional extensibility, the  wrapper reduces the polynomial order of the FLS method, and introduces a sparse multivector which can be further exploited.  The construction and essential properties of the initial involution are explained.

For dimensions 4–7, the optimized inverse requires fewer real multiplications than a same-dimension multivector product.  Dimensions 8–11, follow easily on common hardware with similar efficiency.  For dimensions 12–15, higher-precision arithmetic demonstrates the generality of the algorithm when the arithmetic can support it.  Performance on common hardware is reported.  Directions are included for obtaining a Python package with examples for numerical computation in Clifford algebra, supporting arbitrary signatures including degenerate cases.
\end{abstract}

\maketitle

\section*{Introduction}

Clifford algebra practitioners operating within a few dimensions will be familiar with multivector inverse formulas written in terms of common involutions following graded patterns.  In 2018 Acus and Dargys \cite{acus2018inverse} extended such formulas to 6 dimensional multivectors.  Determinative formulas like these are useful, but bound to a specific dimension.  The following year, Hitzer and Sangwine introduced a method for extending any even dimensional inverse algorithm to the next odd dimension \cite{hitzer2019construction}, thus supporting an inverse method for dimension 7.  Lately, Abdulkhaev and Shirokov \cite{abdulkhaev2025engageSlides} \cite{abdulkhaev2025explicit} revealed an algorithm for dimension 7 that cleverly folds the compound data types described by Hitzer into the involution construction.  This innovation relieves programmers from formal data-type conversions, leaving them grounded in real arrays.

After some vigorous algebraic manipulation of a 6D Acus and Dargys involution formula into polynomial form, this author this year published another 6D multivector inverse algorithm \cite{jones2026another}.  One advantage of the polynomial form was two fewer multivector products required to reach the determinant. The algorithm possessed a beautiful symmetry and simplicity which turned out to have a reason: 
the coefficients are found within a polynomial numerical method called Faddeev--LeVerrier--Souriau (FLS).  The author excitedly reached for an 8D inverse algorithm using the next level of FLS coefficients only to be disappointed.  It turned out, that in order to support the FLS core, the initial involution (reversion), had to be modified.

Working through the details of converting the Abdulkhaev 7D algorithm to polynomial form, the odd step involution innovation that was needed began to show: negate the reversion for the 7D extension.  In higher odd dimensions that formula became, more precisely: append to the block negation pattern its negative. Experimentation uncovered the pattern for the even steps as simply doubling the block size from the previous odd pattern.

With these insights revealed, a general algorithm across dimensions took form.  In the sections to follow the algorithm is explained.  What can be achieved with freely available open-source software, on inexpensive commercial hardware, is documented, and made available online.  

On many occasions, we will reference the pygae (Python Geometric Algebra Enthusiasts) team and their Python package \texttt{clifford}  \cite{python_clifford}.  We will apply the "pygae" label to results obtained from their software.

\section{Algorithm Description}

\subsection{The Envelope}

Most of the involution formulas for the inverse of multivector \(A\) given by Acus and Dargys begin by computing:  \(H = A\widetilde{A}\). 
It follows that \(A^{-1} = \widetilde{A}H^{-1}\) provided \(H\) is invertible. 
An inverse for \(A\) therefore follows from an inverse for \(H\).

\subsection{Equivalent Expressions in Six Dimensions}

Acus and Dargys give the determinant of a 6D multivector \cite{acus2024characteristic} as
\[\det(A) = \frac{1}{3}\left( HH\,\overline{HH} + 2H\,\overline{\overline{H}\;\overline{\overline{H}\;\overline{H}}} \right), \qquad \text{where } \overline{X} = 2\langle X \rangle_0 - X.\]
The bar involution has the effect of negating all but the scalar term in the multivector.  5 multivector products are required to realize this formula.

The polynomial form, written with a scalar-tweak operator, requires only 3 multivector products:
\[\det(A) = -\frac{1}{3}\bigl(((((H\circ 3)H)\circ 1)H)\circ \tfrac{1}{3}\bigr)H, \qquad \text{where } X \circ s = X - (1+s)\langle X\rangle_0 .\]
This is a main result of the author's prior paper \cite{jones2026another}.
The scalar-tweak operator has the effect of multiplying the scalar term of the multivector by $-s.$
Notice that $\overline{X} =- X \circ 1.$ 

Provided that the determinant is different than zero, either of these routes to the determinant offers an expression: 
\[\det(A) = SH = HS, \text{ so that } H^{-1}=\frac{S}{\det(A)}\]
With \(H^{-1}\) we can complete \(A^{-1}=\widetilde{A}H^{-1}\).

\subsection{Modifying Reversion}

Determining an alternative involution to reversion to expand FLS processing to higher dimensions is a main result of the present communication.  The pattern described in Listing \ref{lst:tau} is laid over reversion (by direct product) to determine which multivector components flip sign.  This combination is being called the \(\tau\)-involution, and is a function of dimension.  For dimensions less than seven, reversion is unmodified. Each odd step, beginning at \(d=7\), appends the negative of the previous even step; each subsequent even step doubles the uniform block size of each block of the previous odd step.  

\begin{lstlisting}[float, caption={Signs for the Tau Involution at Dimensions 4-11}, label={lst:tau}, language=Python]
invo_4   = [ 1]*16
invo_5   = [ 1]*32
invo_6   = [ 1]*64
invo_7   = [ 1]*64  + [-1]*64
invo_8   = [ 1]*128 + [-1]*128
invo_9   = [ 1]*128 + [-1]*128 + [-1]*128 + [ 1]*128
invo_10  = [ 1]*256 + [-1]*256 + [-1]*256 + [ 1]*256
invo_11  = ([ 1]*256 + [-1]*256 + [-1]*256 + [ 1]*256
          + [-1]*256 + [ 1]*256 + [ 1]*256 + [-1]*256)
\end{lstlisting}

It is important to realize that the order of components in this implementation is \textit{ordinal}.  In this order the basis-blades have a fixed location, so that the bits of an address correspond to vector-bases.  For instance: the coefficient for element $e_{356}$ will be found at \texttt{Reg[0b110100]} in a translation that treats $e_1$ as the first vector-basis.  Without this statement, the specification of an involution by an array of signs is left ambiguous.  

The block sign pattern is the Thue--Morse sequence: better known in the device of a \textit{parity bit} appended to binary integers, only here converted to a sign under the map: \(\{\oplus:0,1\} \mapsto \{\times:1,-1\}\).  This pattern is perhaps more appropriately described as the \textit{grade involution} pattern which isolates the bivector (even) algebra.  While grade involution is the direct product of reversion and Clifford conjugation, what it means to apply the grade involution only to the upper range of bits, in combination with reversion, is uncertain.  There may be a better way to describe this involution: a subject for further investigation.

As a final remark about the \(\tau\)-involution, it is strange to discover that it is \textbf{not} a graded involution, though its properties are similar to the familiar anti-involutions.

\subsection{Faddeev--LeVerrier--Souriau and its Tweaks}

A general method to reveal the eigen-structure of a multivector can be turned to the inverse problem, as described by Shirokov in theorem 4 of \cite{shirokov2021computing} and implemented in the publicly available Python package \texttt{clifford} \cite{python_clifford}.  This polynomial method, known as Faddeev--LeVerrier--Souriau (FLS), is beautifully simple but suffers from a rapid loss of precision as the order of the polynomial doubles every other additional dimension.  For this reason, the pygae team makes the method available but does not systematically employ it.

The tweak coefficients of the FLS recursion at order \(2^n\) are generally described by the formula:
\[FLS(n) = \left[ \frac{2^n - 1}{1},\ \frac{2^n - 2}{2},\ \ldots,\ \frac{2^{n-1}}{2^{n-1}},\ \ldots,\ \frac{2}{2^n - 2},\ \frac{1}{2^n - 1} \right]\]
which for the first four orders gives
\begingroup
\small
\setlength{\arraycolsep}{2pt}
\[
\begin{array}{r c *{15}{c} c}
FLS(1) = & [ & & & & & & & & 1 & & & & & & & & ] \\[1ex]
FLS(2) = & [ & & & & 3, & & & & 1, & & & & \frac{1}{3} & & & & ] \\[1ex]
FLS(3) = & [ & & 7, & & 3, & & \frac{5}{3}, & & 1, & & \frac{3}{5}, & & \frac{1}{3}, & & \frac{1}{7} & & ] \\[1ex]
FLS(4) = & [ & 15, & 7, & \frac{13}{3}, & 3, & \frac{11}{5}, & \frac{5}{3}, & \frac{9}{7}, & 1, & \frac{7}{9}, & \frac{3}{5}, & \frac{5}{11}, & \frac{1}{3}, & \frac{3}{13}, & \frac{1}{7}, & \frac{1}{15} & ]
\end{array}
\]
\endgroup

While FLS can be directly applied to the source multivector the advantage here is to the reduction in FLS order: half the order for even dimensions and one fourth the order for odd dimensions.  The reduction in order extends the range in dimension before loss of precision foils the practical usefulness of the algorithm.  See Table \ref{tbl:order}.  

\begin{table}[tb]
\centering
\setlength{\tabcolsep}{2pt}
\renewcommand{\arraystretch}{1.05}
\caption{Recursion order, native FLS against \texorpdfstring{$\tau$}{tau}-FLS.}
\label{tbl:order}
\begin{tabular}{lrrrrrrrr}
\toprule
order & 4D & 5D & 6D & 7D & 8D & 9D & 10D & 11D \\
\midrule
Native FLS & \makecell{\(2^2\)\\(4)} & \makecell{\(2^3\)\\(8)} & \makecell{\(2^3\)\\(8)} & \makecell{\(2^4\)\\(16)} & \makecell{\(2^4\)\\(16)} & \makecell{\(2^5\)\\(32)} & \makecell{\(2^5\)\\(32)} & \makecell{\(2^6\)\\(64)} \\
\(\tau\)-FLS & \makecell{\(2^1\)\\(2)} & \makecell{\(2^1\)\\(2)} & \makecell{\(2^2\)\\(4)} & \makecell{\(2^2\)\\(4)} & \makecell{\(2^3\)\\(8)} & \makecell{\(2^3\)\\(8)} & \makecell{\(2^4\)\\(16)} & \makecell{\(2^4\)\\(16)} \\
\bottomrule
\end{tabular}%
\end{table}
Notice that odd dimensions process at the same order as the \textbf{next} even dimension with native FLS, but at the same order as the \textbf{previous} even dimension with \(\tau\)-FLS.  Overall, substituting \(\tau\)-FLS shifts the effort down by three dimensions.

The shift in FLS processing order is an immediate advantage, paid with the cost of two sparse multivector multiplies.  With a little programming effort, exploiting the repetitive sparse multiplications which remain in domain over an FLS sweep realizes a dramatic reduction in effort.

\subsection{The \tautitle-domain}

Certain mathematical properties attributed to reversion are preserved in the \(\tau\)-modification to reversion.  These properties make the optimizations described later possible.

Let \(\tau\) be a multivector involution \(\tau(\tau(A)) = A\), described by component sign changes, which is also an anti-automorphism \(\tau(AB) = \tau(B)\tau(A)\). 
Its \(+1\) fixed set is called the \(\tau\)-domain, and comprises about half the algebra.
Although this domain is not closed under the geometric product, it has four properties that together confine the FLS recursion to it.

\textbf{Entry.} \(C = A\,\tau(A)\) lies in the domain, since \(\tau(C) = \tau(\tau(A))\,\tau(A) = A\,\tau(A) = C\).

\textbf{Powers stay in domain.} \(\tau(X^n) = \tau(X^{n-1})\tau(X) = X^{n-1}X = X^n\).

\textbf{Polynomials stay in domain.} \(\tau\) is linear, so real linear combinations of powers remain in the domain.

\textbf{Commuting products stay in the domain.} For $X, Y$ in the $\tau$-domain, $\tau(XY) = \tau(Y)\,\tau(X) = YX$. Hence $XY$ lies in the $\tau$-domain if and only if $XY = YX$.

\begin{proof}
Let $\mathcal{Z} = \mathbb{R}$ for even $d$ and $\mathcal{Z} = \mathbb{R}\oplus\mathbb{R}I$ for odd $d$. Every element of $\mathcal{Z}$ is central in every signature, including degenerate ones. Degeneracy can enlarge the center, but only by blades that are vacuously central: wherever such a blade fails to commute with a generator, that generator is null and lies in the blade, so both products vanish. These additions play no role here; only $\mathcal{Z}$ is used. Under the $\tau$-construction, $\tau(Z) = Z$ for all $Z \in \mathcal{Z}$, placing $\mathcal{Z}$ entirely within the $\tau$-domain. Each FLS iterate is formed by multiplying by $C$ and shifting by a central element; hence, every iterate is a polynomial in $C$ with coefficients in $\mathcal{Z}$. 

Because central elements commute with all elements of the algebra and powers of $C$ commute with each other ($C^i C^j = C^{i+j} = C^j C^i$), any two such polynomials $X = p(C)$ and $Y = q(C)$ commute: $XY = YX$. Consequently, the commutator product $X \times Y = \tfrac{1}{2}(XY - YX)$ vanishes, and the geometric product reduces to the symmetric product:
\[
XY = \tfrac{1}{2}(XY + YX).
\]
Since $X$ and $Y$ lie in the $\tau$-domain and commute, $\tau(XY) = \tau(Y)\tau(X) = YX = XY$, proving that $XY$ remains in the $\tau$-domain.
\end{proof}

This property, guaranteeing that FLS iterations remain in domain is remarkable given that the \(\tau\)-domain is not in general, closed under the geometric product.  Since \(\tau\)-domain products which stay in domain are commutative, we need record only one sign per qualifying \(\tau\)-domain element pair; and those element pairs are easy to pick out: just those that commute.

\subsection{The Model Algorithm in Python}

To conclude the algorithm description, we attach Listing \ref{lst:model} modeling the algorithm in Python in terms of whole multivectors.  Some explanations are necessary. The multivector \texttt{copy()} method instantiates an independent multivector of the same value.  Otherwise, an aliasing action occurs with assignment \texttt{X=Y}.   \texttt{Reg[]} is a multivector's coefficient array. \texttt{Reg[0]} selects the scalar coefficient, while \texttt{Reg[-1]} selects the last element of the array, the pseudoscalar coefficient.

\begin{lstlisting}[float, caption={Multivector Inverse Model (Python)}, label={lst:model}, language=Python]
Involution[d-4]     # tau = reversion x Thue-Morse layer
Tweak[d-4]          # FLS coefficients, e.g. [3, 1, 1/3] at d = 6,7

def mv_inv(A):                    # 4 <= d <= 11
    d, odd = A.dimensions, A.dimensions & 1

    tA = Involution[d-4] * A      # tA = tau(A)
    C = A * tA                    # C = A tau(A)  -- lands in domain
    
    T = Tweak[d-4]                # FLS sweep, all inside the domain
    D = C.copy()
    D.Reg[0] *= -T[0]
    if odd: D.Reg[-1] *= -T[0]
    for t in T[1:]:
        D *= C
        D.Reg[0] *= -t            # tweak the scalar ...
        if odd: D.Reg[-1] *= -t   # ... and the pseudoscalar
    E = D * C                     # collapses to the center

    if odd:                       # center is E = a + bI
        Ep = E.copy(); Ep.Reg[-1] *= -1 
        rC = (1 / (Ep * E).Reg[0]) * Ep * D
    else:                         # center is a scalar only
        rC = (1 / E.Reg[0]) * D

    return tA * rC                 # 1/A = tau(A) / C
\end{lstlisting}

\section{Optimization}\label{sec:optimize}

The product of two arbitrary d-dimensional multivectors generally requires \(4^d\) real multiplies.  However, no product in this algorithm actually requires that many, when structural particulars are exploited.  Each optimization effectively avoids multiplying by zero, or producing zero when components are known to cancel.

\subsection{Real-Product Counts}\label{subsec:real-prod-counts}

Let's grant that the number of real-products is an objective measure of effort for a numerical procedure.  This gives us an objective comparison between the full-multivector unoptimized model shown in Listing \ref{lst:model}, and an optimized version that reduces the number of real-products to just those required after eliminating multiplication by zero, or to a canceled product known in advance.  It may be that sign-flipping gets translated to a real-product in execution, but we will not count these.  Nor will we count the extra real-products introduced when substituting extended-precision products for double-precision products in dimensions 12-15.  

First we describe the inputs to the real-product count formulas by dimension.  Among these counts are the shapes of the data structures that govern FLS-processing.  "Slots" is the size of the \(\tau\)-domain, and "pairs" is the number of distinct \(\tau\)-input combinations which produce each \(\tau\)-component. The table of pairs \texttt{T[1:s][p]} is of the shape "slots by pairs".  A sign table \texttt{S[1:s][p]} of the same shape indicates the sign associated with each pair.  This sign table bears the effects of choice of signature.  A non-degenerate signature generates a sign table containing \{-1,1\}.  Zeros are introduced to the sign table for degenerate signatures.  You can see the use of these tables in Python Listing \ref{lst.tau_mult}.

 Note that \(c = 1 + (d \bmod 2)\) counts the central tau-slots: the scalar, and at odd $d$ the pseudoscalar. It is independent of signature; vacuously central blades of a degenerate algebra are processed as ordinary slots.
\begin{table}[t]
\centering
\caption{Inputs to the real-product count formulas by dimension \(d\)}
\label{tab:tau-inputs}

\noindent\makebox[\linewidth][c]{%
\begin{minipage}[t]{0.62\linewidth}
\centering
\small
\setlength{\tabcolsep}{3pt}
\renewcommand{\arraystretch}{1.1}
\begin{tabular}{rrrrr}
\toprule
$d$ & $n$ & $c$ & $s$ & $p$ \\
\midrule
4 & 2 & 1 & 6 & 1 \\
5 & 2 & 2 & 12 & 2 \\
6 & 4 & 1 & 28 & 6 \\
7 & 4 & 2 & 56 & 12 \\
\midrule
8 & 8 & 1 & 120 & 28 \\
9 & 8 & 2 & 240 & 56 \\
10 & 16 & 1 & 496 & 120 \\
11 & 16 & 2 & 992 & 240 \\
\midrule
12 & 32 & 1 & 2,016 & 496 \\
13 & 32 & 2 & 4,032 & 992 \\
14 & 64 & 1 & 8,128 & 2,016 \\
15 & 64 & 2 & 16,256 & 4,032 \\
\bottomrule
\end{tabular}
\end{minipage}%
\hspace{0.03\linewidth}%
\begin{minipage}[t]{0.35\linewidth}
\small
\setlength{\tabcolsep}{2pt}
\renewcommand{\arraystretch}{1.1}
\begin{tabular}{@{}ll@{}}
$d$ &= dimensions, \\
$n$ &= FLS order, \\
$c$ &= central slots,\\
$s$ &= slots, \\
$p$ &= pairs.
\end{tabular}
\end{minipage}%
}

\end{table}

We divide the algorithm into three sections: the step-down into the tau-domain, the FLS polynomial iterations, and the step-up from the tau-domain back to a full multivector.  
These are the formulas for the number of required real-products once the zeros in the sparse multivectors have been eliminated.

\begin{align}
M_{\downarrow} &= 2^d + (s-1)\,2^{d-1} = (s+1)\,2^{d-1} \label{eq:step-down}\\
M_{\tau} &= (n-2)\bigl(c + cs + 2p(s-c)\bigr) + cs + 2(c-1) \label{eq:tau-sweep}\\
M_{\uparrow} &= s \cdot 2^d \label{eq:step-up}
\end{align}

\begin{table}[t]
\centering
\caption{Real multiplies per optimized inverse expressed in units of one arbitrary multivector product ($4^d$ real multiplies), along side the model's multivector product count.}
\label{tab:tau-mults}
\begin{tabular}{rrrrrr}
\toprule
$d$ & $M_{\downarrow}$ & $M_{\tau}$ & $M_{\uparrow}$ & total & model \\
\midrule
4 & 0.22 & 0.02 & 0.38 & 0.62 & 3\\
5 & 0.20 & 0.03 & 0.38 & 0.60 & 5 \\
6 & 0.23 & 0.18 & 0.44 & 0.84 & 5\\
7 & 0.22 & 0.18 & 0.44 & 0.84 & 7 \\
\midrule
8 & 0.24 & 0.62 & 0.47 & 1.33 & 9\\
9 & 0.24 & 0.62 & 0.47 & 1.33 & 11 \\
10 & 0.24 & 1.59 & 0.48 & 2.32 & 17\\
11 & 0.24 & 1.59 & 0.48 & 2.32 & 19 \\
\midrule
12 & 0.25 & 3.58 & 0.49 & 4.32 & 33\\
13 & 0.25 & 3.58 & 0.49 & 4.32 & 35 \\
14 & 0.25 & 7.57 & 0.50 & 8.31 & 65\\
15 & 0.25 & 7.57 & 0.50 & 8.31 & 67 \\
\bottomrule
\end{tabular}
\end{table}

\subsection{The Envelope: Entering and Exiting the \tautitle-domain}\label{subsec:envelope}

The \textit{slots} are determined and listed in an array \texttt{I[]}.  This is then the map for stepping into and out of the \(\tau\)-domain.  We are careful to ensure for odd dimensions that the formula \texttt{\char`~I[n]==I[-1-n]} holds.  Here the tilda represents the complement (trimmed to the given dimension); and negative indexing counts from the end of the array.

Because the components of the two multivectors in the product \(C=A\tau(A)\) are redundant, the product of any two distinct components, \(a,b\) of \(A\) will either \(=\pm 2ab\) or \(0\).  Since multiplication by \(\pm 2\) does not require the full machinery of a real multiplication, we don't count it.  We know that whenever the component result is not in domain the product is \(0\) so the products we require are easy to determine.  This is the \textit{step-down cost} described by equation \ref{eq:step-down}.

The \textit{step-up cost} in equation \ref{eq:step-up}, is optimized simply by avoiding the zeros in the \(\tau\)-domain \(C^{-1}\), when computing \(A^{-1}=\tau(A)C^{-1}\)

\subsection{The \tautitle FLS iteration}\label{subsec:tau-fls}

The \(\tau\)-FLS iteration computes the polynomial which reduces to the center of the algebra.  
The multivector products in this section of processing are optimized to a table driven function, described in Python Listing \ref{lst.tau_mult}.  The table \texttt{T[][]} holds an array of \texttt{(a,b)} pairs of indecies into the compressed \(\tau\)-multivector array which produce each non-central component of the product.  The table \texttt{S[][]} holds one sign per pair (or possibly zero in the case of a degenerate signature) and has the same shape as \texttt{T[][]}.  The construction of \texttt{S[][]} encodes the effects of the chosen signature of the algebra.  The signs (or zero) associated with each square of the \(\tau\)-components is tabled in \texttt{S0[]} for computing the scalar part of the \(\tau\)-product.  For odd dimensions a similar array of signs (or zero) associated with the complementary product is tabled in \texttt{SP[]}.

A pair is recorded in the table only when \(e_c=e_ae_b\) stays in the domain, which happens exactly when \(e_a\) and \(e_b\) commute. That is why the two terms of a pair share a single sign. The \(\tau\)-product \(Z=XY\) is then expressed mathematically for even dimensions:

\[Z_c = \sum_{(a,b)\in T_c} S_{c,ab}\bigl(X_aY_b + X_bY_a\bigr), \qquad Z_0 = \sum_{a\in I} S_{0,a}X_aY_a .\]
For odd dimensions, the pseudoscalar component is computed:
\[Z_n = \sum_{a\in I} S_{n,a}X_aY_{\bar a}.\]

Here $n$ is the tau-index of the pseudoscalar and $\bar a$ is the tau-index of the blade complementary to $a$. The central components $Z_0$ and $Z_n$ each take one product per slot, $|I|$ terms, whereas each non-central component $Z_c$ takes two products per pair, $2|T_c|$ terms.

\begin{lstlisting}[float, caption={Tau Multiplication}, label={lst.tau_mult}, language=Python]
# s rows of p index pairs tabled in T[]
# s by p sign table S[][]
# s signs in S0[] for tau-domain squares
def tau_mult( X, Y ): -> Z
    Z[0] = sum(S0[k] * X[k]*Y[k] for k in range(s))    # scalar
    if odd: # pseudoscalar
        Z[-1] = sum(SP[k] * X[k]*Y[-1-k] for k in range(s))
     for r in range(1, s - odd):        # non-central tau-slots
        z = 0.0
        for k, (a, b) in enumerate(T[r]): # one sign, two terms
            z += S[r][k] * (X[a]*Y[b] + X[b]*Y[a])
        Z[r] = z
    return Z
\end{lstlisting}

The last polynomial iteration is extracted from the loop because we want to special case that multivector product knowing that the result is a very sparse multivector containing only the \textit{center}: a lone scalar when the dimension is even, or accompanied by the pseudoscalar when the dimension is odd.  

Referring back to the model in Listing \ref{lst:model}, we now describe optimization of the last FLS product, which is special cased because it collapses to the center:  \(E=DC,\) from the model.  So, the reciprocal we seek is just  \(1/C=D/E\)  for even dimensions where \(E\) is just a scalar.  For odd dimensions, \(E\) includes the pseudoscalar as well, making a little extra processing needed to reduce this to a scalar we can divide by.  For this we invent \(E'\) which is the two component multivector \(E\)  with the pseudoscalar negated.  Then we can compute the scalar value: \(E'E = E[0]*E[0]-I^2*E[-1]*E[-1]\) where \(I^2\in \{-1,0,1\}\) is the pseudoscalar squared.  So, we now have for the odd case: \(E'E=E'DC\), or \(1/C=(E'D)/(E'E)\).  Our real-product count for the last multivector product in the FLS sweep becomes: \(s\) for even dimensions, and \(2+2s\) for odd dimensions (where \(s\) is the "slot count") .  Or by formula: \(cs +2(c-1)\).   Where \(c\) is the number of central slots.

Here is the full formula for the number of real-multiplications required in processing FLS in the \(\tau\)-domain (reiterating equation \ref{eq:tau-sweep}):
\[M_\tau = (n-2)(2p(s-c)+cs+c) + cs+2(c-1)\]
The formula in words:  Each of the \(n-1\) \(\tau\)-products except the last costs two real-products per pair \((2p)\) for each of the non-central slots \((s-c)\), whereas the central slots require one real-product per slot, finishing with one tweak per central slot.  The last \(\tau\)-product produces only the center, so it needs just the \(cs\) real-products landing in the center.  Finally we include the 2 real-product cost to form \(E'E\) when the center is composite.

\subsection{Only Even Dimension Pair Tables are Stored}\label{subsec:even_tbls_only}

For non-degenerate signatures the pseudoscalar lies in the domain and the domain factors:
\[(u_1 + w_1 I)(u_2 + w_2 I) = (u_1u_2 + I^2 w_1w_2) + (u_1w_2 + w_1u_2)I\]
with \(u\) and \(w\) in the even \((d-1)\) domain. The signature of the top generator enters as a single bit, \(I^2 = \pm 1\).
Remember that there is a "complement" organizing device relating the second half of these sign tables to the first half, as expressed in the slot identifier array \texttt{I[]}.

\begin{table}[tb]
\centering
\small
\setlength{\tabcolsep}{4pt}
\renewcommand{\arraystretch}{1.05}
\caption{Stored even-dimension tables and the odd dimensions they derive.}
\label{tbl:storage}
\begin{tabular}{lrrrr}
\toprule
stored table & \(d=4\) & \(d=6\) & \(d=8\) & \(d=10\) \\
\midrule
slots & 6 & 28 & 120 & 496 \\
pairs per row & 1 & 6 & 28 & 120 \\
size & \(<1\) KB & 1 KB & 27 KB & 481 KB \\
derives & \(d=5\) & \(d=7\) & \(d=9\) & \(d=11\) \\
\bottomrule
\end{tabular}
\end{table}

The whole family for \(d = 4 \ldots 11\) occupies 509 KB. 

No odd pair table need ever be built.  However, the sign tables are another matter.  
For degenerate signatures the factorization described above fails: $e_a I = 0$ whenever $e_a$ contains a null generator, so the complement half of the domain cannot be reached through $I$.  
Expanding the sign tables to the full dimensional requirement for degenerate signatures solves this problem.  The full sign table, initialized with appropriate zeros accurately supports degeneracy.

\section{Performance}\label{sec:perf}

The author developed and tested these algorithms in Python 3.12.13 using an Anaconda distribution on an Acer 515 Chromebook Plus. The AMD CPU features \textbf{4 physical cores with SMT enabled (8 logical threads)}, with a base frequency of 2.1 GHz and a dynamic boost frequency up to \textbf{3.7 GHz}. The processor supports \textbf{AVX2 and FMA3 instruction sets}, which are employed via \textbf{Numba's Just-In-Time (JIT) compiler} to optimize high-dimensional multivector operations. The system uses \textbf{8GB of DDR4 RAM}, with \textbf{6.5 GiB} available to the virtualized Linux environment.  The hardware supports a \textbf{64-byte cache line alignment} and \textbf{512 KB of L2 cache per core}.  

As an interpreted language, Python's computing priorities emphasize interactive exploratory development at the expense of computational speed and efficiency.  For numerical work, some of the efficiency of a compiled language can be recovered for numerical processes with declared data types (numpy), and a just-in-time compiler (numba).  The benchmarks here must be understood in this context.
For these reasons, and the intrusion of background processing, benchmarks listed here have to be understood as samplings from the range of expectations.

\subsection{Benchmarks}

Four inverse implementations are bench-marked for comparison.  First the \textbf{model}, as described in abbreviated form in Listing \ref{lst:model}, is timed.  This represents the unoptimized \(\tau\)-FLS algorithm coded in multivector operations. This method transports well to other platforms that support the geometric product.  Next, is the \textbf{optimized \tautitle-FLS} which eliminates induced zeros.  The third method adds a \textbf{Newton-Schulz} refinement to as much as double the precision of the inverse when needed. This at least requires taking the product of the multivector and its computed inverse to assess the residual.  The fourth method could be called \textbf{native FLS} which follows Shirokov's implementation, without first using the \(\tau\)-involution.

Those four methods are supplemented with the independent pygae (Python Geometric Algebra Enthusiasts team) implementation.  The inverse is computed from Hitzer's closed-form formulas for dimensions 4 and 5; and for dimensions without a closed-form formula, pygae relies upon the LAPACK linear algebra package \cite{python_clifford} for an inverse.

Kingdon (v3.0.0) \cite{roelfs2025willingkingdoncliffordalgebra} was also considered. Its built-in inverse did not complete for a general multivector beyond $d = 4$, and FLS written with its product slows sharply from $d = 7$ onward, where Kingdon switches from generated code to direct computation.  We omit Kingdon timings from the tables, concluding that dense high dimensional multivectors are outside its focus.  

For these tests, Euclidean signature was chosen.  Behavior under other signatures may vary.  For instance, mixed or degenerate signatures may increase the possibility of ill-conditioning.  We are considering these issues out of scope for this effort: though baseline regression testing does verify support for arbitrary signature, including degenerate algebras.

\begin{table}[tb]
\centering
\caption{Time / Inverse}
\label{tbl:times}
\small
\setlength{\tabcolsep}{8pt}
\begin{tabular}{rrrrrr}
\toprule
D & Model & $\tau$-FLS & + NS & FLS & pygae \\
\midrule
\multicolumn{6}{l}{\textit{(values in \si{\micro\second})}} \\
4  & 21.17  & 7.135  & 24.21 & 21.27 & 21.52  \\
5  & 42.12  & 8.867  & 29.07 & 53.06 & 56.97 \\
6  & 72.40  & 15.47  & 38.86 & 84.15 & 142.9 \\
7  & 212.1  & 36.94  & 77.18 & 384.3 & 561.2 \\
\midrule
\multicolumn{6}{l}{\textit{(values in \si{\milli\second})}} \\
8  & 0.7262 & 0.1480 & 0.2358 & 1.111 & 4.395 \\
9  & 2.932  & 0.6636 & 1.478  & 7.818 & 26.52 \\
10 & 16.89  & 6.192  & 8.843  & 30.39 & 181.2 \\
11 & 82.98  & 28.14  & 39.08 & 266.8 & 1044 \\
\bottomrule
\end{tabular}
\end{table}

\begin{figure}[tb]
  \centering
  \includegraphics[width=\linewidth]{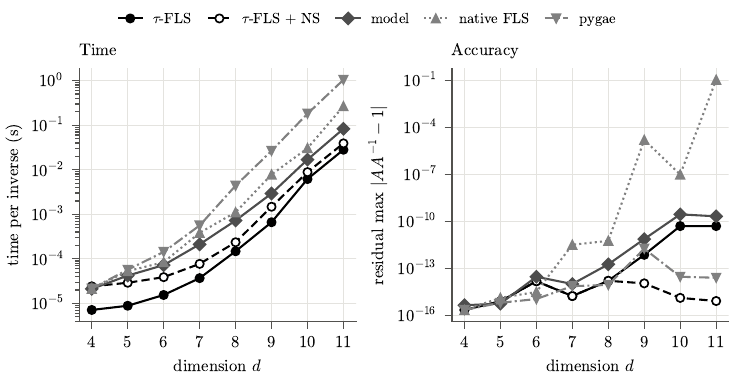}
  \caption{Benchmark Timings and Residuals by Method over Dimensions}
  \label{fig:bench_graph}
\end{figure}

\begin{table}[tb]
\centering
\caption{Residuals -- $\max\lvert A A^{-1} - 1 \rvert$}
\label{tbl:residuals}
\small
\setlength{\tabcolsep}{6pt}
\begin{tabular}{rrrrrr}
\toprule
D & Model & $\uptau$-FLS & + NS & FLS & pygae \\
\midrule
4  & 4.4e-16 & 2.2e-16 & 2.2e-16 & 2.2e-16 & 2.2e-16 \\
5  & 5.6e-16 & 7.8e-16 & 7.8e-16 & 1.3e-15 & 6.3e-16 \\
6  & 2.9e-14 & 1.5e-14 & 1.5e-14 & 3.1e-15 & 1.1e-15 \\
7  & 1.0e-14 & 1.7e-15 & 1.7e-15 & 3.2e-12 & 7.3e-15 \\
8  & 1.8e-13 & 1.6e-14 & 1.6e-14 & 5.7e-12 & 8.7e-15 \\
9  & 7.4e-12 & 7.3e-13 & 1.1e-14 & 1.6e-05 & 1.8e-12 \\
10 & 2.8e-10 & 5.0e-11 & 1.3e-15 & 9.8e-08 & 2.9e-14 \\
11 & 2.1e-10 & 5.0e-11 & 8.2e-16 & 1.1e-01 & 2.5e-14 \\
\bottomrule
\end{tabular}
\end{table}

\begin{table}[tb]
\centering
\caption{optimized $\tau$-FLS improvement factor over other methods}
\label{tbl:perf_ratios}
\small
\setlength{\tabcolsep}{6pt}
\renewcommand{\arraystretch}{1.0}
\begin{tabular}{l*{8}{r}}
\toprule
Method & 4 & 5 & 6 & 7 & 8 & 9 & 10 & 11 \\
\midrule
$\tau$-FLS / model      & 3.0 & 4.8 & 4.7 & 5.7  & 4.9  & 4.4  & 2.7  & 2.9 \\
$\tau$-FLS / native FLS & 3.0 & 6.0 & 5.4 & 10.4 & 7.5  & 11.8 & 4.9  & 9.5 \\
$\tau$-FLS / pygae      & 3.0 & 6.4 & 9.2 & 15.2 & 29.7 & 40.0 & 29.3 & 37.1 \\
\bottomrule
\end{tabular}
\end{table}

The benchmark used a "best of 5" methodology in attempt to avoid the intrusion of background processing. All inverses were pre-warmed to avoid including numba compile time. 

benchmark package versions:
\begin{itemize}
    \item{Python version = 3.12.13} 
    \item{numpy version = 2.4.6} 
    \item{numba version = 0.65.1} 
    \item{clif-ga version = 0.2.0}
    \item{clifford version = 1.5.1}
\end{itemize}

\subsection{Reaching with Extended Precision and Parallel Processing}

With standard double precision losing precision on the \(\tau\)-FLS result at around 12D, 
a version of the inverse in double-double precision was created, and tested in the range 10D-15D.  The results of that test are shown in Table \ref{tab:perf_residual}.
The primary purpose of this extended precision processing was to show that the algorithm is properly general in its mathematical form; and will produce a numerical inverse for a well conditioned multivector provided the arithmetic can support it.
Secondarily, the exercise tests a practical extension to the range in dimension the multivector inverse reaches.

\begin{table}[htbp]
  \centering
  \caption{Performance and residual comparison across dimensions $d$ \\ with FLS orders $N$.}
  \label{tab:perf_residual}

  \setlength{\tabcolsep}{3pt}
  \small
  \resizebox{\textwidth}{!}{%
  \begin{tabular}{rrrrrrrrr}
    \toprule
    $d$ & $N$ &
    \multicolumn{1}{c}{\makecell{Build\\(s)}} &
    \multicolumn{1}{c}{\makecell{f64\\residual}} &
    \multicolumn{1}{c}{\makecell{dbl+dbl\\residual}} &
    \multicolumn{1}{c}{\makecell{f64\\(ms)}} &
    \multicolumn{1}{c}{\makecell{dbl+dbl\\(ms)}} &
    \multicolumn{1}{c}{Ratio} &
    \multicolumn{1}{c}{\makecell{Peak\\(MB)}} \\
    \midrule
    10 & 16 & 0.019 & $4.84 \times 10^{-10}$ & $8.53 \times 10^{-15}$ &   6.694 &    19.777 & 2.955 & 241 \\
    11 & 16 & 0.012 & $2.04 \times 10^{-11}$ & $7.61 \times 10^{-15}$ &  28.374 &    79.153 & 2.790 & 371 \\
    12 & 32 & 0.251 & $2.32 \times 10^{-04}$ & $1.75 \times 10^{-15}$ & 166.381 &   637.633 & 3.832 & 371 \\
    13 & 32 & 0.257 & $2.95 \times 10^{-04}$ & $7.18 \times 10^{-15}$ & 695.906 &  2551.612 & 3.667 & 371 \\
    14 & 64 & 6.648 & \multicolumn{1}{c}{--} & $8.42 \times 10^{-08}$ & \multicolumn{1}{c}{--} & 20837.411 & \multicolumn{1}{c}{--} & 589 \\
    15 & 64 & 7.523 & \multicolumn{1}{c}{--} & $4.39 \times 10^{-07}$ & \multicolumn{1}{c}{--} & 92922.446 & \multicolumn{1}{c}{--} & 590 \\
    \bottomrule
  \end{tabular}%
  }
\end{table}

The build times for constructing the required tables gets to be noticeable at higher dimensions, but is small in comparison to the wait times for extended precision inverse results (7 seconds vs 1.5 minutes at 15D).  (Realize that build times are one time costs whereas inverse times are per-inverse costs.)

With long wait times for an extended precision result, the potential speed up using parallel processing becomes important.  Numba was coached to use parallel processing for both the double and the double-double versions of FLS processing; and returned the timings reported in Table \ref{tab:perf_residuals_parallel}.

\begin{table}[t]
  \centering
  \caption{Performance and residual comparison after enabling numba parallel execution.}
  \label{tab:perf_residuals_parallel}

  \setlength{\tabcolsep}{3pt}
  \small
  \resizebox{\textwidth}{!}{%
  \begin{tabular}{rrrrrrrrr}
    \toprule
    $d$ & $N$ &
    \multicolumn{1}{c}{\makecell{Build\\(s)}} &
    \multicolumn{1}{c}{\makecell{f64\\residual}} &
    \multicolumn{1}{c}{\makecell{dbl+dbl\\residual}} &
    \multicolumn{1}{c}{\makecell{f64\\(ms)}} &
    \multicolumn{1}{c}{\makecell{dbl+dbl\\(ms)}} &
    \multicolumn{1}{c}{Ratio} &
    \multicolumn{1}{c}{\makecell{Peak\\(MB)}} \\
    \midrule
    10 & 16 & 0.019 & $4.84 \times 10^{-10}$ & $8.53 \times 10^{-15}$ &   1.910 &     4.800 & 2.513 & 252 \\
    11 & 16 & 0.016 & $2.04 \times 10^{-11}$ & $7.61 \times 10^{-15}$ &   7.789 &    19.250 & 2.471 & 379 \\
    12 & 32 & 0.253 & $2.32 \times 10^{-04}$ & $1.75 \times 10^{-15}$ &  52.567 &   157.498 & 2.996 & 379 \\
    13 & 32 & 0.261 & $2.95 \times 10^{-04}$ & $7.18 \times 10^{-15}$ & 216.595 &   690.291 & 3.187 & 379 \\
    14 & 64 & 6.722 & \multicolumn{1}{c}{--} & $8.42 \times 10^{-08}$ & \multicolumn{1}{c}{--} &  6605.000 & \multicolumn{1}{c}{--} & 598 \\
    15 & 64 & 6.554 & \multicolumn{1}{c}{--} & $4.39 \times 10^{-07}$ & \multicolumn{1}{c}{--} & 36246.962 & \multicolumn{1}{c}{--} & 598 \\
    \bottomrule
  \end{tabular}%
  }
\end{table}

The just-in-time numba compiler provides an easy parallel processing option for today's multi-core processors (8 threads over 4 cores on the author's chromebook).  The comparative results for parallel processing in the two numeric precisions is provided in Table \ref{tab:parallel_speedup}.

\begin{table}[t]
  \centering
  \caption{Speedup from enabling numba parallel execution (relative to non-parallel results).}
  \label{tab:parallel_speedup}
  \begin{tabular}{rrr}
    \toprule
    $d$ & \multicolumn{1}{c}{f64 speedup} & \multicolumn{1}{c}{dbl+dbl speedup} \\
    \midrule
    10 & 3.50 & 4.12 \\
    11 & 3.64 & 4.11 \\
    12 & 3.17 & 4.05 \\
    13 & 3.21 & 3.70 \\
    14 & \multicolumn{1}{c}{--} & 3.15 \\
    15 & \multicolumn{1}{c}{--} & 2.56 \\
    \bottomrule
  \end{tabular}
\end{table}

Engaging numba's parallel processing saves enough time to warrant its use whether standard double precision or extended precision is employed.  The current baseline software uses extended precision for dimensions 12D-15D, standard double precision for dimensions 4D-11D, and parallel processing for the whole range of optimized inverse processing, 4D-15D.  Please realize that earlier timings in this paper are for single-threaded processes only.  If needed, parallel processing is easily switched off.

\section{Conclusion}

The author is pleased to have demonstrated a practical multivector inverse through 15D for use on common hardware with open-source software.
\begin{itemize}
    \item 0D-3D: heuristic formulas
    \item 4D-7D: improving level of effort over previously known formulas
    \item 8D-11D: extending capability without resorting to matrices
    \item 12D-15D: practical processing with extended precision and parallel processing
\end{itemize}

Aside from the invention of another way to compute the inverse of a multivector, there is a certain satisfaction in uncovering the construction of the \(\tau\)-involutions.  The feeling that a deeper source of mathematical vitality has been struck comes through in the slots and pairs counts in Table \ref{tab:tau-inputs}.  Notice that the number of pairs at any dimension equals the number of slots two dimensions back.  \(p[d+2]=s[d].\)  This is an unintentional property that derives from an experimental construction.  This observation suggests that we should find \(\tau\)-involutions at the lesser dimensions such that: \(s[2]=1,\) and \(s[3]=2;\) and, there \textbf{is} one: \textit{Clifford conjugation} as described by Lounesto \cite{lounesto2001clifford}, pg.57.  With conjugation as our \(\tau\)-involution at these dimensions, \(A\tau(A)\) is exactly the center.  Since we routinely remove the center from the table of pairs, we conclude that there are no pairs at \(d=2,3.\)  We conclude further that there are neither slots nor pairs at dimensions 0 or 1 where the multivector is already central.  As shown in Table \ref{tab:tau-pairless}, the cases for dimension < 4 logically belong above Table \ref{tab:tau-inputs}.  However, since these dimensions show no information useful to FLS processing this table annex was held back until now.

\begin{table}[t]
\centering
\caption{Inputs to the real-product counts formulas for the pair-less dimensions, and the associated \(\tau\)-involution}
\label{tab:tau-pairless}

\noindent\makebox[\linewidth][c]{%
\begin{minipage}[t]{0.62\linewidth}
\centering
\small
\setlength{\tabcolsep}{3pt}
\renewcommand{\arraystretch}{1.1}
\begin{tabular}{rrrrrc}
\toprule
$d$ & $n$ & $c$ & $s$ & $p$ & \(\tau\)-involution \\
\midrule
0 & 0 & 1 & 0 & 0 & \texttt{+} \\
1 & 0 & 2 & 0 & 0 & \texttt{+ +} \\
2 & 0 & 1 & 1 & 0 & \texttt{+ - - -} \\
3 & 0 & 2 & 2 & 0 & \texttt{+ - - - - - - +} \\
\bottomrule
\end{tabular}
\end{minipage}%
\hspace{0.03\linewidth}%
\begin{minipage}[t]{0.35\linewidth}
\small
\setlength{\tabcolsep}{2pt}
\renewcommand{\arraystretch}{1.1}
\begin{tabular}{@{}ll@{}}
$d$ &= dimensions, \\
$n$ &= FLS order, \\
$c$ &= central slots,\\
$s$ &= slots, \\
$p$ &= pairs.
\end{tabular}
\end{minipage}%
}

\end{table}

To conclude with an overall assessment of this package's performance, we offer the following benchmarks in Table \ref{tbl:shipped}.
\begin{table}[tb]
\centering
\caption{As-shipped timing: one geometric product and one inverse, $\mathrm{Cl}(d,0)$, random multivectors.}
\label{tbl:shipped}
\small
\setlength{\tabcolsep}{8pt}
\begin{tabular}{rrrrr}
\toprule
$d$ & product & inverse & \makecell{inverse /\\product} & residual \\
\midrule
\multicolumn{5}{l}{\textit{(values in \si{\nano\second})}} \\
 0 & 1848 & 1792 & 0.97 & 0 \\
 1 & 1894 & 1966 & 1.04 & $1.1 \times 10^{-16}$ \\
 2 & 1920 & 2305 & 1.20 & $2.8 \times 10^{-17}$ \\
 3 & 2096 & 2368 & 1.13 & $5.7 \times 10^{-17}$ \\
\midrule
\multicolumn{5}{l}{\textit{(values in \si{\micro\second})}} \\
 4 & 2.337 & 6.173 & 2.64 & $2.2 \times 10^{-16}$ \\
 5 & 3.615 & 8.054 & 2.23 & $7.8 \times 10^{-16}$ \\
 6 & 7.521 & 14.52 & 1.93 & $1.5 \times 10^{-14}$ \\
 7 & 20.23 & 35.81 & 1.77 & $1.7 \times 10^{-15}$ \\
 8 & 63.90 & 150.3 & 2.35 & $1.6 \times 10^{-14}$ \\
 9 & 226.9 & 664.4 & 2.93 & $7.3 \times 10^{-13}$ \\
\midrule
\multicolumn{5}{l}{\textit{(values in \si{\milli\second})}} \\
10 & 0.862 & 3.302 & 3.83 & $5.0 \times 10^{-11}$ \\
11 & 3.458 & 14.77 & 4.27 & $5.0 \times 10^{-11}$ \\
12 & 14.54 & 172.9 & 11.89 & $1.7 \times 10^{-14}$ \\
13 & 148.9 & 1159 & 7.78 & $1.3 \times 10^{-14}$ \\
14 & 239.5 & 6375 & 26.62 & $1.9 \times 10^{-9}$ \\
15 & 1040 & 28524 & 27.42 & $1.0 \times 10^{-6}$ \\
\bottomrule
\end{tabular}
\end{table}

The assumption that a real product count was a good measure of effort does not bear out with these timings.  But, they are honest, and they represent a sampling of what you can expect from this \cite{clifga_pypi_0_2_0} package if you choose to use it.

\section*{Statements and Declarations}

\subsection*{Funding}
The author received no financial support for the research, authorship, and/or publication of this article.

\subsection*{Competing Interests}
The author declares no competing interests.

\subsection*{Data Availability}
The author developed and tested these algorithms using an Anaconda distribution of Python and Jupyter Notebook. Original code and documentation are available from the public repository (see Ref.~\cite{jones_clif_2026}).  The reader will find browser friendly PDF versions of the Jupyter notebook files there.

The software package "clif-ga" can be imported from PyPI \cite{clifga_pypi_0_2_0}.  Documentation accompanying the software can be found at ReadTheDocs \cite{clifga_docs}.  This is a capable package with special focus on the geometric product and its inverse.  Aside from the inverse support described in this paper, the geometric product is implemented with virtual sign-tables culled from Hadamard matrices that leave memory for the application.  The author would entertain collaborators on this project.

\FloatBarrier

\bibliographystyle{unsrturl}
\bibliography{clifford}

\end{document}